\documentclass[12pt]{article}
\usepackage{amssymb,amsthm,amsmath,amsfonts,latexsym,tikz,hyperref,ytableau}
\usepackage[hmargin=1in,vmargin=1in]{geometry}

\newtheorem{thm}{Theorem}[section]
\newtheorem{prop}[thm]{Proposition}
\newtheorem{cor}[thm]{Corollary}
\newtheorem{lem}[thm]{Lemma}
\newtheorem{conj}[thm]{Conjecture}
\newtheorem{exa}[thm]{Example}

\newcommand{\fI}{\mathfrak I}
\DeclareMathOperator{\RS}{RS}
\DeclareMathOperator{\hook}{hook}
\DeclareMathOperator{\dhook}{dhook}
\DeclareMathOperator{\2row}{2row}
\DeclareMathOperator{\d2row}{d2row}
\DeclareMathOperator{\skm}{skm}
\DeclareMathOperator{\ecol}{ecol}

\newcommand{\ben}{\begin{enumerate}}
\newcommand{\een}{\end{enumerate}}
\newcommand{\ble}{\begin{lem}}
\newcommand{\ele}{\end{lem}}
\newcommand{\bth}{\begin{thm}}
\renewcommand{\eth}{\end{thm}}
\newcommand{\bpr}{\begin{prop}}
\newcommand{\epr}{\end{prop}}
\newcommand{\bco}{\begin{cor}}
\newcommand{\eco}{\end{cor}}
\newcommand{\bcon}{\begin{conj}}
\newcommand{\econ}{\end{conj}}
\newcommand{\bde}{\begin{defn}}
\newcommand{\ede}{\end{defn}}
\newcommand{\bex}{\begin{exa}}
\newcommand{\eex}{\end{exa}}
\newcommand{\barr}{\begin{array}}
\newcommand{\earr}{\end{array}}
\newcommand{\btab}{\begin{tabular}}
\newcommand{\etab}{\end{tabular}}
\newcommand{\beq}{\begin{equation}}
\newcommand{\eeq}{\end{equation}}
\newcommand{\bea}{\begin{eqnarray*}}
\newcommand{\eea}{\end{eqnarray*}}
\newcommand{\bal}{\begin{align*}}
\newcommand{\bce}{\begin{center}}
\newcommand{\ece}{\end{center}}
\newcommand{\bpi}{\begin{picture}}
\newcommand{\epi}{\end{picture}}
\newcommand{\bpp}{\begin{picture}}
\newcommand{\epp}{\end{picture}}
\newcommand{\bfi}{\begin{figure} \begin{center}}
\newcommand{\efi}{\end{center} \end{figure}}
\newcommand{\bprf}{\begin{proof}}
\newcommand{\eprf}{\end{proof}\medskip}
\newcommand{\capt}{\caption}
\newcommand{\bsl}{\begin{slide}{}}
\newcommand{\esl}{\end{slide}}
\newcommand{\bfr}{\begin{frame}}
\newcommand{\efr}{\end{frame}}

\newcommand{\hqed}{\hfill \qed}

\newcommand{\hso}[1]{\hspace{-1pt}}

\newcommand{\qmq}[1]{\quad\mbox{#1}\quad}

\newcommand{\sbe}{\subseteq}

\newcommand{\ptn}{\vdash}

\def\<{\langle}
\def\>{\rangle}

\newcommand{\ree}[1]{(\ref{#1})}

\newcommand{\ra}{\rightarrow}

\newcommand{\io}{\iota}

\newcommand{\la}{\lambda}

\newcommand{\si}{\sigma}

\newcommand{\La}{\Lambda}

\newcommand{\ba}{{\bf a}}
\newcommand{\bb}{{\bf b}}

\newcommand{\bbR}{{\mathbb R}}

\newcommand{\bbZ}{{\mathbb Z}}

\newcommand{\cL}{{\cal L}}

\newcommand{\fS}{{\mathfrak S}}

\DeclareMathOperator{\id}{id}

\DeclareMathOperator{\sh}{sh}

\DeclareMathOperator{\SYT}{SYT}

\begin{document}
\pagestyle{plain}

\title{Longest increasing subsequences and log concavity
}
\author{
Mikl\'os B\'ona\\[-5pt]
\small Department of Mathematics, University of Florida,\\[-5pt]
\small P.O. Box 118105, Gainesville, FL 32611-8105, USA, {\tt bona@ufl.edu}\\
Marie-Louise Lackner\thanks{Supported by the Austrian Science Foundation FWF, grant P25337-N23}\\[-5pt]
\small Institute of Discrete Mathematics and Geometry, Technische Universit\"at Wien\\[-5pt]
\small Wiedner Hauptstra\ss e 8-10, 1040 Vienna, Austria, {\tt  marie-louise.lackner@tuwien.ac.at}\\
Bruce E. Sagan\\[-5pt]
\small Department of Mathematics, Michigan State University,\\[-5pt]
\small East Lansing, MI 48824-1027, USA, {\tt sagan@math.msu.edu}
}

\date{\today\\[10pt]
	\begin{flushleft}
	\small Key Words: hook shape, involution, log concavity, longest increasing subsequence, permutation, Robinson-Schendsted algorithm, skew-merged permutation, standard Young tableau, Tracy-Widom distribution, two-rowed shape
	                                       \\[5pt]
	\small AMS subject classification (2010):  05A05
	\end{flushleft}}

\maketitle

\begin{abstract}
Let $\pi$ be a permutation of $[n]=\{1,\dots,n\}$ and denote by $\ell(\pi)$  the length of a longest increasing subsequence of $\pi$.  Let $\ell_{n,k}$ be the number of permutations $\pi$ of $[n]$ with $\ell(\pi)=k$.  Chen conjectured that the sequence $\ell_{n,1},\ell_{n,2},\dots,\ell_{n,n}$ is log concave for every fixed positive integer $n$.  We conjecture that the same is true if one is restricted to considering involutions and we show that these two conjectures are closely related.  We also prove various analogues of these conjectures concerning permutations whose output tableaux under the Robinson-Schensted algorithm have certain shapes.  In addition, we present a proof of Deift that part of the limiting distribution is log concave.  Various other conjectures are discussed.
\end{abstract}

%
%

\section{Introduction}

Let $\fS_n$ be the symmetric group of all permutations of $[n]=\{1,2,\dots,n\}$.  We will view $\pi=\pi_1 \pi_2 \dots \pi_n\in\fS_n$ as a sequence (one-line notation).  Let $\ell(\pi)$ denote the length of a longest increasing subsequence of $\pi$.  For example, if $\pi=4172536$ then $\ell(\pi)=4$ because  the subsequence $1256$ is increasing and there is no longer such sequence.  Define
$$
L_{n,k}=\{\pi\in\fS_n\ :\ \ell(\pi)=k\} \qmq{and}  \ell_{n,k}=\#L_{n,k}
$$
where the hash symbol denotes cardinality.  

The statistic $\ell(\pi)$ plays an important role in a number of combinatorial contexts, for example in famous theorems of Erd\H{o}s and Szekeres~\cite{es:cpg} and of Schensted~\cite{sch:lid}.  
The problem of determining the distribution of $\ell(\pi)$ in a random permutation $\pi$ of length $n$ was solved in a tour de force by Baik, Deift, and Johansson~\cite{bdj:dli}. 
The history of this problem is described by Aldous and Diaconis~\cite{aldous}; see also the recently published book by Romik \cite{romik} on the subject. 

The statistic $\ell(\pi)$ is not only interesting from a combinatorial or algorithmic point of view:  It is also connected with biology via the Ulam distance which is used to model evolutionary distance in DNA research~\cite{ula:ipb}.  The {\em Ulam distance} between $\pi,\si\in\fS_n$, denote $U(\pi,\si)$, is the minimum number of steps needed to obtain $\si$ from $\pi$ where a step consists of taking an element of a sequence and placing it somewhere else in the sequence.  If $\id$ is the identity permutation then it is easy to see that $U(\id,\pi)=n-\ell(\pi)$.  Indeed, one can fix a longest increasing subsequence of $\pi$ and then move all the elements of $\id$ which are not in that subsequence to the appropriate places to form $\pi$.  

A sequence of real numbers $l_1,l_2,\dots,l_n$ is said to be {\em log concave} if $l_{k-1} l_{k+1}\le l_k^2$ for all $k\in [n]$.  Here we use the convention that $l_0=l_{n+1}=0$.  Log concave sequences appear often in algebra, combinatorics, and geometry; see the survey articles of Stanley~\cite{sta:lus} and Brenti~\cite{bre:lus}.  It is interesting to note that there are many other ways to define distance in molecular biology, and these sequences are typically either known to 
be log concave, or conjectured to be log concave. See the book \cite{genome} for a collection of examples.

Our main object of study is the following conjecture which appeared in an unpublished manuscript of William Chen from 2008.
\bcon[\cite{che:lqc}]
\label{ell_nk}
For any fixed $n$, the sequence
$$
\ell_{n,1}, \ell_{n,2},\dots,\ell_{n,n}
$$
is log concave.
\econ
We have verified this conjecture for $n\le 50$ by computer and will give other evidence for its truth below.

Let $\fI_n$ denote the set of involutions in $\fS_n$, i.e., those permutations whose square is the identity.  Also define
$$
I_{n,k}=\{\pi\in\fI_n\ :\ \ell(\pi)=k\} \qmq{and}  i_{n,k}=\#I_{n,k}.
$$
We conjecture the following.
\bcon
\label{i_nk}
For any fixed $n$, the sequence
$$
i_{n,1}, i_{n,2},\dots,i_{n,n}
$$
is log concave.
\econ
Again, this conjecture has been verified for $n\le 50$.  We will show that these two conjectures are closely related.

The rest of this paper is structured as follows.  In the next section we will see that there is a close connection between Conjectures~\ref{ell_nk} and~\ref{i_nk}.  We will also derive another relation between sequences counting certain permutations and those counting certain involutions.  Section~\ref{htr} restricts attention to permutations whose output tableaux under the Robinson-Schensted map have certain shapes.    Fixed-point free involutions are considered in Section~\ref{fpf}.  Baik, Deift, and Johansson~\cite{bdj:dli} proved that, with suitable scaling, the sequence in Conjecture~\ref{ell_nk} converges to the Tracy-Widom distribution as $n\ra\infty$.  In Section~\ref{twd}, we present a proof of Deift that this distribution is log concave for nonnegative $x$ where $x$ is the independent variable.  We end with more conjectures related to~\ref{ell_nk} and~\ref{i_nk}.

\section{Involutions}
\label{inv}

To connect Conjectures~\ref{ell_nk} and~\ref{i_nk}, we will need some properties of the Robinson-Schensted correspondence.  For more information about this important map, see the the texts of Sagan~\cite{sag:sym} or Stanley~\cite{sta:ec2}.

Let $\la=(\la_1,\la_2,\dots,\la_k)$ be a partition of $n$, written $\la\ptn n$.  
We denote by $\SYT \la$ the set of all standard Young tableaux of shape $\la$, and if $P\in\SYT\la$, then we will also write $\sh P=\la$.
The Robinson-Schensted map is a bijection
$$
\RS:\fS_n \ra \bigcup_{\la\ptn n} (\SYT\la)^2,
$$
i.e., a permutation of length $n$ is identified with a pair of standard Young tableaux of size $n$ and of the same shape.
We will use the notation $\RS(\pi)=(P,Q)$.  Since $\sh P =\sh Q$ we can define the {\em shape of $\pi$} to be the common shape of its output tableaux.  It will also be convenient to define for  pairs of permutations $\sh(\pi,\si)=(\sh\pi,\sh\si)$.

We need two important results about $\RS$, the first due to Schensted~\cite{sch:lid} and the second to Sch\"utzenberger~\cite{sch:rcs}.
\bth
\label{RS}
The map $\RS$ has the following properties.
\ben
\item[(1)] If $\sh\pi=(\la_1,\dots,\la_k)$ then $\ell(\pi)=\la_1$.  Also, the length of a longest decreasing subsequence of $\pi$ is the number of cells in the first column of $\la$.
\item[(2)] A permutation $\pi$ is an involution if and only if $\RS(\pi)=(P,P)$ for some standard Young tableau $P$.\hqed
\een
\eth

By (2), there is a canonical bijection between involutions and standard Young tableaux.  Because of this, we will go freely back and forth between involutions and tableaux without further mention.

One way to prove Conjecture~\ref{ell_nk} would be to find injections $F:L_{n,k-1}\times L_{n,k+1}\ra L_{n,k}^2$ for all $n,k$.  Call any map with this domain and range {\em shape preserving} if
$$
\sh(\pi,\pi')=\sh(\si,\si') \implies \sh F(\pi,\pi') = \sh F(\si,\si')
$$
for all $(\pi,\pi'), (\si,\si')\in L_{n,k-1}\times L_{n,k+1}$.  We will also apply this terminology to functions 
$f:I_{n,k-1}\times I_{n,k+1}\ra I_{n,k}^2$.  Our first result shows that if one can prove Conjecture~\ref{i_nk} using a shape-preserving injection, then one gets Conjecture~\ref{ell_nk} for free.
\bth
\label{sp}
Suppose that there is a shape-preserving injection $f:I_{n,k-1}\times I_{n,k+1}\ra I_{n,k}^2$ for some $n,k$.  Then there is a shape-preserving injection $F:L_{n,k-1}\times L_{n,k+1}\ra L_{n,k}^2$.
\eth
\bprf
Given $f$, one can construct $F$ as the composition of the following maps.
\begin{align*}
(\pi,\pi')&\stackrel{RS^2}{\longmapsto} \left( (P,Q),\ (P',Q')\right)\\[5pt]
&\longmapsto \left( (P,P'),\ (Q,Q')\right)\\
&\stackrel{f^2}{\longmapsto} \left(( S,S'),\ (T,T')\right)\\[5pt]
&\longmapsto \left( (S,T),\ (S',T')\right)\\
& \hspace{-0.2cm} \stackrel{(RS^{-1})^2}{\longmapsto} \left(\si,\si'\right)
\end{align*}
Note that in applying $f^2$ we are treating each of the tableaux $P,\dots,T'$ as involutions in the manner discussed earlier.  Also, the fact that $f$ is shape preserving guarantees that $\sh S =\sh T$ and $\sh S' = \sh T'$ so that one can apply the inverse Robinson-Schensted map at the last stage.
\eprf

There is another way to relate the log concavity of sequences such as those in Conjectures~\ref{ell_nk} and~\ref{i_nk}.  Let $\La$ be a set of partitions of $n$.  Define
$$
L_{n,k}^\La = \{\pi\in L_{n,k}\ |\ \sh\pi\in\La\} \qmq{and} \ell_{n,k}^\La=\# L_{n,k}^\La.
$$
Similarly define $I_{n,k}^\La$ and $i_{n,k}^\La$.  Clearly our original sequences are obtained by choosing $\La$ to be all partitions of $n$.  At the other extreme, there is also a nice relationship between the log concavity of these two sequences.
\ble
\label{La}
Let $\La$ contain at most one partition with first row of length $k$ for all $1\le k\le n$.  Then $\ell_{n,1}^\La,\dots,\ell_{n,n}^\La$ is log concave if and only if $i_{n.k}^\La,\dots,i_{n,k}^\La$ is log concave.
\ele
\bprf
The hypothesis on $\La$ and Theorem~\ref{RS} imply that $\ell_{n,k}^\La=(i_{n,k}^\La)^2$.  The result now follows since the square function is increasing on nonnegative values.
\eprf

\section{Hooks and two-rowed tableaux}
\label{htr}

If one considers permutations whose output tableaux under $\RS$ have a certain shape, it becomes easier to prove log-concavity results analogous to Conjectures~\ref{ell_nk} and~\ref{i_nk}.  Using the notation developed at the end of the previous section, let $\La=\hook$ be the set of all partitions of $n$ that have the shape of a hook. That is, these are Ferrers shape of a partition of $n$ 
in which the first part is $k$, and there are $n-k$ additional parts, each equal to 1,  where $1\le k\le n$.
These shapes will be denoted by $(k, 1^{n-k})$.

\bth
\label{hook}
For any fixed $n$, the sequences
$$
\ell_{n,1}^{\hook}, \ell_{n,2}^{\hook},\dots,\ell_{n,n}^{\hook}
\qmq{and}
i_{n,1}^{\hook}, i_{n,2}^{\hook},\dots, i_{n,n}^{\hook}
$$
are log concave.
\eth

\bfi
\begin{tikzpicture}
\draw[gray] (0,0) grid (5,5);
\draw[ultra thick,dashed] (1,1)--(3,1)--(3,2)--(4,2)--(4,3);
\node at (.55,.7) {$(a,b)$};
\end{tikzpicture}
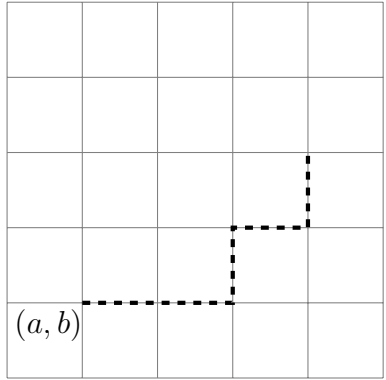
\capt{ The path $p=EENEN$   \label{p}}
\efi

\bprf
Since the set ``$\hook$" satisfies the hypothesis of Lemma~\ref{La}, it suffices to prove the involution result.  For an algebraic proof note that, by the hook formula or direct counting,
$$
i_{n,k}^{\hook}=\mbox{number of partitions of shape $(k,1^{n-k})$}= \binom{n-1}{k-1}.
$$
It is well known and easy to prove by cancellation of factorials that this sequence of binomial coefficients is log concave.

There is also a standard combinatorial proof of the log concavity of this sequence using the technique of Lindstr\"om~\cite{lin:vri}, later used to great effect  by Gessel and Viennot~\cite{gv:bdp}.  We review it here for use in the proof of the next theorem. A {\em NE-lattice path}, $p$, starts at a point $(a,b)\in\bbZ^2$ and takes unit steps north and east, denoted $N$ and $E$ respectively.  See Figure~\ref{p} for an illustration.  There is a bijection between the set of partitions $P$ with 
$\sh P=(k,1^{n-k})$ and NE-lattice paths from $(a,b)$ to $(a+k-1,b+n-k)$ where the $i$-th step of $p$ is $E$ if and only if $i+1$ is in the first row of $P$.  Our example path corresponds to the tableau
$$
P=\begin{ytableau}
1&2&3&5\\
4\\
6
\end{ytableau}\ .
$$
To construct an injection $f:I_{n,k-1}^{\hook}\times I_{n,k+1}^{\hook}\ra (I_{n,k}^{\hook})^2$ we interpret a pair of involutions in the domain as a pair of lattice paths $(p,p')$ where $p$ goes from $(1,0)$ to $(k-1,n-k+1)$ and $p'$ goes from $(0,1)$ to $(k,n-k)$.  The paths $p$ and $p'$ must intersect, since $p$ starts on the southwest side of $p'$ and ends on the northeast side of it.  Let $z$ be the last (most northeast) point in which $p$ and $p'$ intersect.  We then map this pair to $(q,q')$ where $q$ follows $p$ up to $z$ and then follows $p'$, and vice-versa for $q'$.  It is now a simple matter to show that this gives a well-defined injection $f$.
\eprf

\begin{figure}
    \centering
    \hspace{0.02\textwidth}
    \begin{minipage}{.45\textwidth}
        \centering
        \begin{tikzpicture}[scale=1]

\draw[dashed] (1,-1)--( 2,-1)--(3,-1)--( 3,0.05)--(4,0.05)--(4,1);
\foreach \x/\y in {1/-1, 2/-1, 3/-1, 3/0, 4/0,4/1}
\node[draw=black,fill=white,rectangle, inner sep=1.5mm] at (\x,\y) {};

\foreach \x/\y in {0/0, 1/0, 2/0, 3/0, 4/0, 5/0}
\node[fill=black,circle, inner sep=0.8mm] at (\x,\y) {};

\draw (0,0)--(1,0)--(2,0)--(3,0) -- (4,0)--(5,0);

\node[anchor=north] at (1.5,-1) {$p$};
\node[anchor=north] at (0.5,0) {$p'$};
;

\node[anchor=north west] at (4,0) {$z$};
    
\draw[dotted] (0,0)--(2,2);
\draw[loosely dotted] (1.2,-0.8)--(4,2);
\end{tikzpicture}
    \end{minipage}%
    \hfill
    \begin{minipage}{0.45\textwidth}
        \centering
     \begin{tikzpicture}[scale=1]

\draw[dashed] (1,-1)--( 2,-1)--(3,-1)--( 3,0.05)--(4,0.05)--(4,0)--(5,0);
\foreach \x/\y in {1/-1, 2/-1, 3/-1, 3/0, 4/0,5/0}
\node[draw=black,fill=white,rectangle, inner sep=1.5mm] at (\x,\y) {};

\foreach \x/\y in {0/0, 1/0, 2/0, 3/0, 4/0, 4/1}
\node[fill=black,circle, inner sep=0.8mm] at (\x,\y) {};

\draw (0,0)--(1,0)--(2,0)--(3,0) -- (4,0)--(4,1);

\node[anchor=north] at (1.5,-1) {$q$};
\node[anchor=north] at (0.5,0) {$q'$};
;

\node[anchor=north west] at (4,0) {$z$};
    
\draw[dotted] (0,0)--(2,2);
\draw[loosely dotted] (1.2,-0.8)--(4,2);
\end{tikzpicture}
    \end{minipage}%
    \hfill
    \begin{minipage}{0.45\textwidth}
        \centering
     \begin{tikzpicture}[scale=0.8]
\end{tikzpicture}   
    \end{minipage}
    \hspace{0.05\textwidth}
\caption{The main step in the lattice path proof of Theorem~\ref{thm:2row} }
\label{fig:2row_paths}
\end{figure}
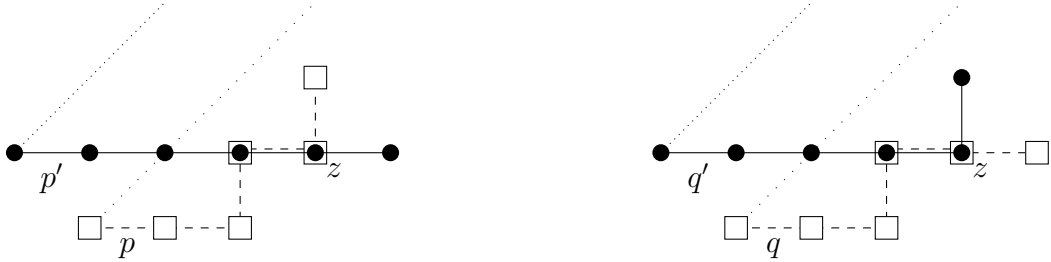

We next turn our attention to shapes with at most two rows.  Let $\La=\2row$ be the set of  shapes of the form $(k,n-k)$ as $k$ varies.  Note that 
$ \ell_{n,k}^{\2row}=0$ for $k<n/2$, but we will still start our sequences at $k=1$ for simplicity.
\bth
For any fixed $n$, the sequences
$$
\ell_{n,1}^{\2row}, \ell_{n,2}^{\2row},\dots,\ell_{n,n}^{\2row}
\qmq{and}
i_{n,1}^{\2row}, i_{n,2}^{\2row},\dots, i_{n,n}^{\2row}
$$
are log concave.
\label{thm:2row}
\eth
\bprf
The  arguments used to prove Theorem~\ref{hook} can be used here as well.  One only needs to be careful about the lattice path proof.  First of all, one maps a tableau to a lattice path using all the elements of the first row (including $1$) for the $E$ steps, and those of the second row for the $N$ steps.  Returning to our example path in Figure~\ref{p}, the corresponding tableau is now
$$
P=\begin{ytableau}
1&2&4\\
3&5
\end{ytableau}\ .
$$
The lattice paths which correspond to $2$-rowed tableaux are the Dyck paths, those which never go above the line of slope $1$ passing through the initial point $(a,b)$.  
See Figure~\ref{fig:2row_paths} for an illustration.
One must now check that if $(p,p')$ maps to $(q,q')$ that $q$ and $q'$ are still Dyck.  This is clear for $q'$ since the portion of $p$ which it uses lies below the line $y=x-1$ and so certainly lies  below $y=x+1$.  For $q$ one must use the fact that $z$ is the last point of intersection.  Indeed, it is easy to see that if $q'$ is not Dyck then there would have to be an intersection point of $p$ and $p'$ later than $z$ which is a contradiction.
\eprf

We remark that the previous result is related to pattern avoidance.
If $\pi\in\fS_k$ is a permutation called the {\em pattern} then a permutation $\si$ {\em avoids} $\pi$ if there is no subsequence $\si'$ of $\si$ of length $k$ which standardizes to $\pi$ when one replaces the smallest element of $\si'$ by $1$, the next smallest by $2$, and so forth.
So, by Theorem~\ref{RS} (1), the  permutations whose output tableaux under RS consist of at most two rows are exactly those that avoid the pattern $321$. 
Thus Theorem~\ref{thm:2row} can be expressed in terms of $321$-avoidance.

Let us now turn to a class of permutations closely related to those of hook shape.  A permutation is {\em skew merged} if it is a shuffle of an increasing permutation and a decreasing permutation.   These permutations have been characterized by Stankova~\cite{sta:fs}  and  the generating fuction for this class was found by Albert and Vatter~\cite{av:ge3}. Again using part (1) of Theorem~\ref{RS}, we see that if $\pi$ is skew merged then $\sh\pi$ has first row of length at least $k$ and first column of length at least $n-k$ for some $k$.  It follows that $\sh\pi$ is either a hook or the union of a hook and the box in the second row and column.  On the other hand, if $\sh\pi$ is a hook then similar reasoning shows that $\pi$ must be skew merged.  However,  not all permutations whose shape is of the second type are skew merged.  For example, both $2413$ and $2143$ have shape $(2,2)$ but the first one is skew merged while the second is not.  We will use a superscript $\skm$ in our notation to restrict to the set of skew merged permutations.  
\bcon
The sequence
$$
\ell_{n,1}^{\skm}, \ell_{n,2}^{\skm},\dots,\ell_{n,n}^{\skm}
$$
is log concave.
\econ
This conjecture has been verified for $n\le 9$.  The reason for the difference in the size of this bound and the previous ones is that for the earlier conjectures we were able to use the hook formula to speed up computations considerably.

It is natural to ask about the analogue of the previous conjecture for involutions.  But it turns out that we have already answered this question in Theorem~\ref{hook}.  To see why, we need the following result of Stankova.  
\bth[\cite{sta:fs}]
\label{sta}
A permutation is skew merged if and only if it avoids the patterns $2143$ and $3412$.\hqed
\eth

\bfi
\begin{tikzpicture}
\draw (0,0)--(4,0)--(4,4)--(0,4)--(0,0) (2,0)--(2,4) (0,2)--(4,2);
\filldraw (2,4) circle(.1);
\filldraw (4,2) circle(.1);
\draw (2,4.4) node{$(j,n)$};
\draw (4.7,2) node {$(n,j)$};
\draw (1,1) node {$A$};
\draw (3,1) node {$B$};
\draw (1,3) node {$C$};
\draw (3,3) node {$D$};
\end{tikzpicture}
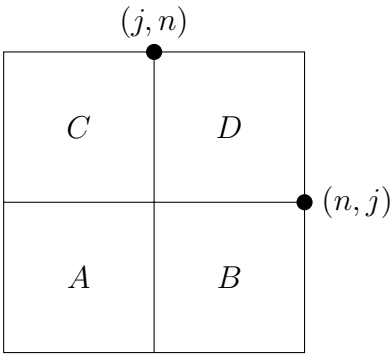
\capt{The diagram of $\io$ \label{io}}
\efi

\bco
Let $\io$ be an involution.  Then $\io$ is skew merged if and only if $\io$ is of hook shape.
\eco
\bprf
We have already observed that the backwards direction holds for all permutations.  For the forward implication, we induct on $n$ where $\io\in\fS_n$. There are two cases depending on whether $n$ is a fixed point or is in a two cycle of $\io$.  In the first case, we have the concatenation $\io=\io' n$.  Thus $\sh\io$ is just $\sh\io'$ with a box for $n$ appended at the end of the first row, and we are done.

In the second case, suppose $n$ is in a cycle with $j<n$.  We represent $\io=\io_1\dots \io_n$ as the set of points $(i,\io_i)$, $i\in[n]$, in the first quadrant of the plane.  The vertical line through $(j,n)$ and the horizontal line through $(n,j)$ divide the box containing $\io$ into four open areas as displayed in Figure~\ref{io}.
Since $\io$ avoids $2143$ by Theorem~\ref{sta}, the points in area $A$ must be increasing.  Since $\io$ is an involution, if there is a point in area $B$ then there must be a corresponding point in area $C$ and vice-versa.  However,  this would contradict the fact that $\io$ also avoids $3412$.  So areas $B$ and $C$ are empty.  This implies that area $A$ contains exactly the fixed points $1,2,\dots,j-1$.  By induction, the involution in area $D$ has shape $(a,1^b)$ for some $a,b$.   It is now easy to see that $\sh\io=(a+j-1,1^{b+2})$ finishing the proof.
\eprf

\section{Fixed-point free involutions}
\label{fpf}

Chen's manuscript also included some conjectures about perfect matchings.  These correspond in $\fS_n$ to involutions without fixed points.  Their output shapes are characterized by the following result of Sch\"utzenberger.

\bth[\cite{sch:cr}]
If $\io$  is an involution then the number of fixed points of $\io$ equals the number of columns of odd length of $\sh\io$.\hqed
\eth

Let $\La=\ecol$ be the set of partitions of $n$ all of whose column lengths are even.  Chen's perfect matching conjecture can now be stated as follows.  The parity of the columns forces the number of elements permuted to be even.
\bcon[\cite{che:lqc}]
For any fixed $n$, the sequence
$$
i_{2n,1}^{\ecol}, i_{2n,2}^{\ecol},\dots,i_{2n,2n}^{\ecol}
$$
is log concave.
\econ
This conjecture has been verified for $2n\le 80$.  Given the development so far, it is natural to also make the following conjecture.
\bcon
For any fixed $n$, the sequence
$$
\ell_{2n,1}^{\ecol}, \ell_{2n,2}^{\ecol},\dots,\ell_{2n,2n}^{\ecol}
$$
is log concave.
\econ
A computer has tested this conjecture for $2n\le 80$.

For certain subsets of $\ecol$ it is possible to prove log concavity results.  We define the {\em double} of a partition $\la=(\la_1,\la_2,\dots,\la_l)$ to be $\la^2=(\la_1,\la_1,\la_2,\la_2,\dots,\la_l,\la_l)$.  First consider the set $\La=\dhook$ of double hooks, namely the partitions of shape $(k^2,1^{2n-2k})$ where $1\le k\le n$.
\bth
\label{dhook}
For any fixed $n$, the sequences
$$
\ell_{2n,1}^{\dhook}, \ell_{2n,2}^{\dhook},\dots,\ell_{2n,n}^{\dhook}
\qmq{and}
i_{2n,1}^{\dhook}, i_{2n,2}^{\dhook},\dots, i_{2n,n}^{\dhook}
$$
are log concave.
\eth
\bprf
Appealing to Lemma~\ref{La} again, we only need to prove the statement about involutions.  Applying the hook formula gives the following expression in terms of a multinomial coefficient
$$
i_{2n,k}^{\dhook}=\frac{1}{(2n-k)(2n-k+1)}\binom{2n}{1,k-1,k,2n-2k}.
$$
Substituting this into the defining inequality for log concavity and cancelling shows that it suffices to prove
\begin{align*}
& (k-1)(2n-2k)(2n-2k-1)(2n-k+1)(2n-k)\\
& \qquad \le(k+1)(2n-2k+2)(2n-2k+1)(2n-k+2)(2n-k-1).
\end{align*}
Substituting $d=2n-2k$ one can write both sides of this inequality as a polynomial in $k$ and $d$.  Moving all terms to the right-hand side we see that one must show
\begin{align*}
2 d^4 + 8 d^3 k + 10 d^2 k^2 + 4 d k^3 + 4 d^3 + 10 d^2 k + 10 d k^2 + 2 k^3 + 2 d^2 + 2dk + 4 k^2 - 4d - 2k -4 \ge0.
\end{align*}
Since $d\ge0$ and $k\ge1$ it is easy to see that the last inequality is valid.
\eprf

We next turn to two-rowed partitions which have been doubled.  Let $\La=\d2row$ be partitions of shape $(k^2,(n-k)^2)$ as $k$ varies.
\bth
For any fixed $n$, the sequences
$$
\ell_{2n,1}^{\d2row}, \ell_{2n,2}^{\d2row},\dots,\ell_{2n,n}^{\d2row}
\qmq{and}
i_{2n,1}^{\d2row}, i_{2n,2}^{\d2row},\dots, i_{2n,n}^{\d2row}
$$
are log concave.
\eth
\bprf
As usual, it suffices to consider the case concerning  involutions.  The hook formula gives, for $n/2\le k\le n$,
$$
i_{2n,k}^{\d2row}=\frac{(2k-n+1) (2k-n+2)^2 (2k-n+3)}{ (k+1)^2 (k+2)^2 (k+3) (n-k+1)}\binom{2n}{k,k,n-k,n-k}.
$$
One now follows the proof of the previous result, except that one uses the substitutions $d=2k-n$ and $c=k-n$.  The result is a polynomial in $c,d$ which has only positive coefficients and so we are done.
\eprf

\section{The Tracy-Widom distribution}
\label{twd}

We now investigate the limiting distribution of the sequence in Conjecture~\ref{ell_nk}.  The {\em Tracy-Widom distribution}, first studied by Tracy and Widom \cite{tracy1994level}, has cumulative distribution function
\beq
\label{F(x)}
F(x)=\exp\left(-\int_x^\infty (t-x) u^2(t) dt\right)
\eeq
where $u(x)$ is a solution to the Painlev\'e II equation
\beq
\label{PII}
u''(x)= x u(x) + 2 u^3(x).
\eeq
satisfying
\beq
\label{bc}
u(x), u'(x) \ra 0 \qmq{as} x\ra\infty.
\eeq
Call any twice-differentiable function  $f:\bbR\ra\bbR$ {\em concave at $x$} if $f''(x)\le 0$.  Similarly, if $D\sbe \bbR$ then we say $f$ is {\em concave on $D$} if it is concave at any point of $D$.  If $f$ only takes on positive values then we will say the function is {\em log concave} at $x$ if the function $\log f$ is concave at $x$.  If $f$ is log concave on an interval of radius one about $x$ then one can prove that $f(x-1) f(x+1)\le f(x)^2$.

Percy Deift (personal communication) has shown that  the density function of the  Tracy-Widom distribution is log concave for  nonnegative $x$.  We thank him for his kind permission to reproduce his proof here.
\bth[Deift]
\label{TW}
The density function of the Tracy-Widom distribution is log concave for $x\ge0$.
\eth
\bprf
To get the  density function for the Tracy-Widom distribution, one must take the derivative $F'(x)$ where $F(x)$ is as given in equation~\eqref{F(x)}.  Then to determine log concavity, we need to compute $(\log F'(x))''$.  Sraightforward calculations give
$$
(\log F'(x))'' = - \frac{u^2(x) h^2(x) + 2 u(x) u'(x) h(x) + u^4(x)}{h^2(x)}
$$
where 
\beq
\label{h(x)}
h(x) = \int_x^\infty u^2(t) dt.
\eeq
So it suffices to show that
\beq
\label{uh}
u^2(x) h^2(x) + 2 u(x) u'(x) h(x) + u^4(x)\ge 0.
\eeq
Multiplying the Painlev\'e II equation~\ree{PII} by $u'(x)$ and integrating from $x$ to infinity one obtains, with the help of the boundary conditions~\ree{bc},
\beq
\label{u'2}
(u'(x))^2 = x u^2(x) + h(x) + u^4(x).
\eeq

We now make the following claims.
\ben
\item[(i)]  We have $u'(x) u(x) < 0$ for $x\ge0$.
\item[(ii)]  We have $h(x)\le u^2(x)$ for $x\ge 1/4$.
\een

To prove the first claim we first assume, towards a contradiction, that  $u'(x)=0$ for some $x\ge0$.  But then equation~\eqref{u'2} forces $u(x)=0$ as well.  It follows that $u(x)$ must be the zero function since it is a solution to a second order equation which is our desired contradiction.

Now, since $u'(x)$ is continuous, we must have $u'(x)>0$ for all $x\ge0$ or $u'(x)<0$ for all $x\ge0$.  Without loss of generality, we can assume the latter since the substitution of $-u(x)$ for $u(x)$ takes a solution of~\ree{PII} into another solution.  We now break the proof of (i) into two cases, the first being when $u(x)\neq 0$ for the given value of $x$.  Assume, towards a contradiction, that $u'(x) u(x) > 0$.  This forces $u(x)<0$.  But then $u(t)<u(x)<0$ for $t>x$ because we have $u'(t)<0$.  This contradicts the boundary condition in~\ree{bc} that $u(t)\ra 0$ as $t\ra\infty$.  So this case can not occur.
The other case is when $u(x)=0$.  But since $u'(x)<0$ there is a $y>x$ where $u(y)<0$ and now we continue as in the previous case.  This finishes the proof of (i).  Note that in the process we have shown that $u(x)>0$ for $x\ge0$.

To prove (ii), it suffices to show that the function $f(x)=u^2(x)-h(x) $ is nonnegative.  Using the definition of $h(x)$, we have
$$
f'(x) = 2 u(x) u'(x) + u^2(x) = u(x) (2 u'(x) + u(x)).
$$
Using equation~\ree{u'2} we see that $(u'(x))^2 > x u^2(x) \ge u^2(x)/4$ for $x\ge1/4$.  Recalling that $u'(x)<0$, this translates to $2u'(x)+u(x)<0$.  Combining this with $u(x)>0$, we see that $f'(x)<0$.  Also the boundary conditions and definition of $h(x)$ show that $f(x)\ra 0$ as $x\ra \infty$.   Thus $f(x)\ge 0$ as desired.

We are now ready to prove the theorem itself.  Dividing equation~\eqref{uh} by $u^2(x)$, it suffices to prove that
$$
g(x) \stackrel{\rm def}{=} h^2(x) - 2 v(x) h(x) + u^2(x) > 0
$$
where
$$
v(x)\stackrel{\rm def}{=} - u'(x)/u(x) >0
$$
by Claim (i).
We now compute, omitting the independent variable and using the definitions of $h(x)$ and $v(x)$,
\begin{align*}
g'&= -2hu^2 - 2v'h +2v u^2 + 2uu'\\
&=-2hu^2 +2 h \frac{u''u-(u')^2}{u^2} -2u'u+2u u'\\
&=-2h\frac{u^4-u''u+(u')^2}{u^2}.
\end{align*}
Using first the Painlev\'e II equation~\ree{PII} and then equation~\eqref{u'2} on the numerator gives
$$
u^4-u''u+(u')^2 = u^4 - x u^2 - 2 u^4 + (u')^2= h.
$$
Thus $g'(x) = -2h^2(x) / u^2(x) <0$.

We claim that $g(x)\ra 0$ as $x\ra \infty$ which will finish the proof since, together with $g'(x)<0$ for $x\ge0$, this forces $g(x)> 0$.  We first need to analyze what happens to $v(x)$.  By equation~\eqref{u'2} we have $v^2 = x + u^2 + h/u^2$.  As $x\ra\infty$ we know that $u^2(x)\ra 0$ and,  by Claim (ii), $h(x)/u^2(x)\le 1$.  Thus $v(x)\sim \sqrt{x}$ and $-2v(x) h(x) \ra 0$ since it is known that $h(x)$ decreases as $\exp(-x^{3/2})$.
We already know that the other two summands in $g(x)$ approach zero, so we are done with the claim and with the proof of the theorem.
\eprf

\section{Other conjectures}
\label{oc}
 
\subsection{More on the Tracy-Widom distribution}

Approximate values for the mean and variance of the  Tracy-Widom distribution are  $\mu\approx -1.77$ and $\si^2\approx 0.81$.  So Theorem~\ref{TW}
only addresses what is going on in the tail.  It would be very interesting to see what can be said  for negative $x$.  Note that the point in the proof where the bound $x\ge0$ was imposed is in the proof of Claim (i). (Claim (ii) is only needed to analyze what goes on as $x\ra\infty$ and so that bound is not the controlling one.)  In particular, we need $x$ to be nonnegative when using equation~\ree{u'2} to conclude that $u'(x)\neq 0$.

If log concavity of the whole distribution could be proved, then one might be able to prove Conjecture~\ref{ell_nk} as follows.  First one would try to find a bound $N$ such that for  $n\ge N$ the log concavity of the Tracy-Widom distribution implies the log concavity of the $\ell_{n,k}$ sequence.  Then, if $N$ were not too large, one could check log concavity for $n<N$ by computer.

\subsection{Real zeros}

Let $a_0,a_1,\dots,a_n$ be a sequence of real numbers.  Consider the corresponding generating function $a(q)=a_0+a_1 q+\dots + a_n q^n$.  It is well known that if the $a_k$ are positive then $a(q)$ having only real roots implies that the original sequence is log concave.  So one might ask about the polynomials $\ell_n(q)$ and $i_n(q)$ for the sequences  in Conjectures~\ref{ell_nk} and~\ref{i_nk}, respectively.  (The fact that they both begin with a zero just contributes a factor of $q$.) Unfortunately, neither seem to be real rooted in general.  In particular, this is true of $\ell_{12}(q)$ and $i_4(q)$.

\subsection{$q$-log convexity}

The $\ell_n(q)$ discussed in the previous subsection do seem to enjoy another interesting property.  We partially order polynomials with real coefficients by letting $f(q)\le g(q)$ if  $g(q)-f(q)$ has nonnegative coefficients.  Equivalently, for any power $q^i$ its coefficient in $f(q)$ is less than or equal to its coefficient in $g(q)$.  Define a sequence $f_1(q), f_2(q), \dots$ to be {\em $q$-log convex} if $f_{n-1}(q) f_{n+1}(q) \ge f_n^2(q)$ for all $n\ge2$.  Another conjecture from the paper of Chen is as follows.
\bcon[\cite{che:lqc}]
The sequence $\ell_1(q),\ell_2(q),\dots$ is $q$-log convex.
\econ
We have verified this conjecture up to $n=50$.  Interestingly  the corresponding conjecture for the involution sequence is false as 
$i_3(q) i_5(q) \not\ge i_4^2(q)$.

\subsection{Infinite log concavity}

There is another way in which one can generalize log concavity.  The {\em $\cL$-operator} applied to a sequence $\ba=(a_0,a_1, \dots, a_n)$ returns a sequence $\bb=(b_0,b_1,\dots,b_n)$ where $b_k=a_k^2-a_{k-1} a_{k+1}$ for $1\le k \le n$ (with the convention that $a_{-1}=a_{n+1}=0$).  Clearly $\ba$ being log concave is equivalent to $\bb=\cL(\ba)$ being nonnegative.  But now one can apply the operator multiple times.  Call $\ba$ {\em infinitely log concave} if $\cL^i(\ba)$ is a nonnegative sequence for all $i\ge1$.  This concept was introduced by Boros and Moll~\cite{bm:ii}  and has since been studied by a number of authors including Br\"and\'en~\cite{bra:isg}, McNamara and Sagan~\cite{ms:ilc}, and Uminsky and Yeats~\cite{uy:uri}.   Here is yet another conjecture of Chen.
\bcon[\cite{che:lqc}]
For any fixed $n$, the sequence
$$
\ell_{n,1}, \ell_{n,2},\dots,\ell_{n,n}
$$
is infinitely log concave.
\econ

Note that proving that a given sequence is infinitely log concave can not necessarily be done by computer since one needs to apply the $\cL$-operator an infinite number of times.  However, we are able to prove the previous conjecture for small $n$ by using a technique developed by McNamara and Sagan~\cite{ms:ilc}.  Given a real number $r$ we say that the sequence $\ba$ is {\em $r$-factor log concave} if $a_k^2\ge r a_{k-1} a_{k+1}$ for all $k$.  So the sequence will be log concave as long as $r\ge1$.  Let $r_0=(3+\sqrt{5})/2$.
\ble[\cite{ms:ilc}]
Let $\ba$ be a sequence of nonnegative real numbers.  If $\ba$ is $r_0$-factor log concave, then so is $\cL(\ba)$.   If follows that if $\ba$ is $r_0$-factor log concave then it is also infinitely log concave. \hqed
\ele

Using this lemma, we can try to prove that $\ba$ is infinitely log concave as follows.  Apply the $\cL$-operator a finite number of times, checking that at each stage the sequence is nonnegative.  If the sequence becomes $r_0$-factor log concave after the finite number of applications, then the lemma allows us to stop and conclude infinite log concavity.   Using this technique and a computer we have proved the following.
\bpr
 The sequence
$$
\ell_{n,1}, \ell_{n,2},\dots,\ell_{n,n}
$$
is infinitely log concave for $n\le50$.
\epr

Finally, one could wonder about infinite log concavity of the involution sequence.  But again, it behaves differently and is not infinitely log concave starting at $n=4$.

\newcommand{\etalchar}[1]{$^{#1}$}

\end{document}